\documentclass[11pt]{amsart}
\usepackage{amsfonts,amsmath,amsthm,amscd,amssymb,latexsym,cite,verbatim,texdraw,floatflt,
caption2}
\RequirePackage{hyperref}
\usepackage[a4paper]{geometry}

\title{Linearly  ordered coarse spaces  }
\author{ Igor  Protasov}

\address{I.Protasov: Taras Shevchenko National University of Kyiv, Department of Computer Science and Cybernetics, Academic Glushkov pr. 4d, 03680 Kyiv, Ukraine}
\email{i.v.protasov@gmail.com}

\begin{document}
\begin{abstract} 
A coarse space 
$X$, endowed with a linear order
 compatible with the coarse structure of $X$,
 is called linearly  ordered.
 We prove that every linearly  ordered coarse space 
 $X$ is locally convex  
  and the asymptotic dimension of $X$  is either $0$ or $1$. If $X$ is metrizable then the family of all right bounded subsets of $X$ has a selector.
\end{abstract}
\maketitle



\section{ Introduction and preliminaries}

Let $\mathcal{K}$ be a class of coarse spaces. 
Given  $X\in \mathcal{K}$, how can one detect whether there exists a linear order on $X$, compatible  with  the coarse structure of $X$?
We used selectors to answer this question if 
$\mathcal{K}$ is one of the following classes: discrete  coarse spaces  \cite{b6}, \cite{b7};
finitary  coarse spaces of groups 
\cite{b8}; 
finitary  coarse spaces of graphs 
\cite{b9}.

In this paper, we continue the investigations of the structure of a linearly ordered coarse space initiated in  \cite{b7}.

In Section 2, we  prove that every linearly ordered coarse space is locally convex, but the coarse structure of $X$ needs not  to be interval. Given a linear order $\leq$ on a set $X$, we characterize the minimal and  maximal coarse structures on $X$, compatible with the interval  bornology of $(X, \leq)$.

In Section 3, we prove that the asymptotic dimension of a linearly ordered coarse space is either 0 or 1.

In Section 4, we construct a selector of the family of right bounded subsets of a metrizable lineary ordered coarse space.

We conclude the paper with Section 5  of comments and open questions.

\vspace{5 mm}

We recall some basic definitions. 
Given a set $X$, a family $\mathcal{E}$  of subsets of $X\times X$ is called a
{\it  coarse structure} on $X$ if

\begin{itemize}
\item{} each $E \in \mathcal{E}$  contains the diagonal $\bigtriangleup _{X}:=\{(x,x): x\in X\}$ of $X$;
\vspace{3 mm}

\item{}  if  $E$, $E^{\prime} \in \mathcal{E}$  then  $E \circ E^{\prime} \in \mathcal{E}$  and
$ E^{-1} \in \mathcal{E}$,    where  $E \circ E^{\prime} = \{  (x,y): \exists z\;\; ((x,z) \in E,  \ (z, y)\in E^{\prime})\}$,    $ E^{-1} = \{ (y,x):  (x,y) \in E \}$;
\vspace{3 mm}

\item{} if $E \in \mathcal{E}$ and  $\bigtriangleup_{X}\subseteq E^{\prime}\subseteq E$  then  $E^{\prime} \in \mathcal{E}$.
\end{itemize}

Elements $E\in\mathcal E$ of the coarse structure are called {\em entourages} on $X$.

For $x\in X$  and $E\in \mathcal{E}$, the set $E[x]:= \{ y \in X: (x,y)\in\mathcal{E}\}$ is called the {\it ball of radius  $E$  centered at $x$}.
Since $E=\bigcup_{x\in X}( \{x\}\times E[x]) $, the entourage $E$ is uniquely determined by  the family of balls $\{ E[x]: x\in X\}$.
A subfamily ${\mathcal E} ^\prime \subseteq\mathcal E$ is called a {\em base} of the coarse structure $\mathcal E$ if each set $E\in\mathcal E$ is contained in some $E^\prime \in\mathcal E^\prime$.

The pair $(X, \mathcal{E})$  is called a {\it coarse space}  \cite{b13} or  a {\em ballean} \cite{b10}, \cite{b11}. 

%We note that  coarse spaces can be considered as asymptotic counterparts of uniform spaces, see [10, Section 1.1].

A coarse  spaces
$(X, \mathcal{E})$ is  called 
{\it connected} if, for any $x, y \in X$, there exists $E\in \mathcal{E}$ such that $y\in E[x]$. 
A subset  $Y\subseteq  X$  is called {\it bounded} if $Y\subseteq E[x]$ for some $E\in \mathcal{E}$,
  and $x\in X$.
If  $(X, \mathcal{E})$ 
is connected  then 
the family $\mathcal{B}_{X}$ of all bounded subsets of $X$  is a bornology on $X$.
We recall that a family $\mathcal{B}$  of subsets of a set $X$ is a {\it bornology}
if $\mathcal{B}$ contains the family $[X] ^{<\omega} $  of all finite subsets of $X$
 and $\mathcal{B}$  is closed   under finite unions and taking subsets. A bornology $\mathcal B$ on a set $X$ is called {\em unbounded} if $X\notin\mathcal B$.
A subfamily  $\mathcal B^{\prime}$ of $\mathcal B$ is called a base for $\mathcal B$ if, for each $B \in \mathcal B$, there exists $B^{\prime} \in \mathcal B^{\prime}$ such that $B\subseteq B^{\prime}$.

Each subset $Y\subseteq X$ defines a {\it subspace}  $(Y, \mathcal{E}|_{Y})$  of $(X, \mathcal{E})$,
 where $\mathcal{E}|_{Y}= \{ E \cap (Y\times Y): E \in \mathcal{E}\}$.
A  subspace $(Y, \mathcal{E}|_{Y})$  is called  {\it large} if there exists $E\in \mathcal{E}$
 such that $X= E[Y]$, where $E[Y]=\bigcup _{y\in Y} E[y]$.

Let $(X, \mathcal{E})$, $(X^{\prime}, \mathcal{E}^{\prime})$
 be  coarse spaces. 
 A mapping $f: X \to X^{\prime}$ is called
  {\it  macro-uniform }  if for every $E\in \mathcal{E}$ there
  exists $E^{\prime}\in \mathcal{E}^{\prime}$  such that $f(E(x))\subseteq  E^{\prime}(f(x))$
    for each $x\in X$.
If $f$ is a bijection such that $f$  and $f ^{-1 }$ are macro-uniform, then   $f  $  is called an {\it asymorphism}.
If  $(X, \mathcal{E})$ and  $(X^{\prime}, \mathcal{E}^{\prime})$  contain large  asymorphic  subspaces, then they are called {\it coarsely equivalent.}

For a coarse space 
$(X,\mathcal{E})$, 
we denote by 
$exp \ X$
 the family  of all non-empty  subsets of $X$
  and by $exp \ \mathcal{E}$ the coarse structure on $exp \ X$
  with the base $\{ exp \ E : E\in  \mathcal{E}\}$, where
 $$(A,B)\in exp \ E 
  \Leftrightarrow A \subseteq E[B], \ \ B\subseteq E[A],$$
and say that $(exp \ X, exp \ \mathcal{E} )$ is the {\it hyperballean} of 
$(X,\mathcal{E})$. 

%For hyperballeans, see \cite{b4}, \cite{b10}, \cite{b11}.

Let  $\mathcal{F}$ be 
  a non-empty  subspace  of $exp \ X$. 
 We say that a macro-uniform mapping 
 $f: \mathcal{F} \longrightarrow X$ 
 is an $\mathcal{F}$-{\it selector} 
 of $(X,\mathcal{E})$ if $f(A)\in A$ for each $A\in \mathcal{F}$. In the case $\mathcal{F}\in [X]^2$,
 $\mathcal{F}= \mathcal{B}_X$
 and $\mathcal{F}= exp \  X$, an $\mathcal{F}$- selector is called a $2$-{\it selector}, 
 a {\it bornologous selector}  
 and a {\it global selector}  respectively.

  We recall that a connected coarse space $(X,\mathcal{E})$ is {\it discrete} if, for each $E\in  \mathcal{E}$, there exists a bounded subset $B$ of $(X,\mathcal{E})$ such that $E[x]=\{x\}$  for each $x\in X\setminus B$. Every bornology $\mathcal{B}$ on a set $X$ defines the discrete coarse space $X_\mathcal{B} = (X,\mathcal{E} _\mathcal{B})$, where  $\mathcal{E}_\mathcal{B}$ is a coarse structure with the base 
$\{ E_B: B\in  \mathcal{B}\}$, $E_B [x]=B$ if $x\in B$ and 
$E_B [x]= \{x\}$
if $x\in X\setminus B$. On the other hand,  every discrete coarse space $(X, \mathcal{E})$ coincides with 
$X_\mathcal{B}$, where $\mathcal{B}$ is the bornology of bounded subsets of $(X, \mathcal{E})$.

%\vspace{3 mm}

\section{ Local convexity and interval bases }

Let $(X, \mathcal{E})$ be a coarse space.
Following \cite{b7},  we say that a linear order 
$\leq$ or $X$ is {\it compatible} with the coarse structure $\mathcal{E}$ if one of the following equivalent conditions holds

\vspace{7 mm}

\begin{itemize}
\item{} for every $E\in \mathcal{E}$, there exists $F\in \mathcal{E}$ such that if $x<y$ and $y\in X\setminus F[x]$ then $x^\prime < y$ for each $x^\prime \in E[x]$ ;

\vspace{3 mm}

\item{}  for every $E\in \mathcal{E}$, there exists $H\in \mathcal{E}$ such that if 
$y<x$ and $y\in X\setminus 
H[x]$ then 
$y< x^\prime $ for each $x^\prime \in E[x]$ ;

\vspace{3 mm}

\item{}  for every $E\in \mathcal{E}$, there exists $K\in \mathcal{E}$ such that if 
$x<y$ and $y\in X\setminus 
K[x]$ then $x^\prime < y^\prime $ for all $x^\prime \in E[x]$,  $y^\prime \in E[y]$.

\end{itemize}

\vspace{5 mm}

A coarse space $(X, \mathcal{E})$,
 endowed with a linear order $\leq$ 
 compatible with $\mathcal{E}$ is called 
 {\it linearly ordered}. 
 In this case, by [7, Proposition 2], the  mapping 
 $f: [X]^2 \longrightarrow X$, defined by 
 $f(A)=min \ A$ is a 2-selector of $(X, \mathcal{E})$
 and if  $(X, \mathcal{E})$ is connected then each interval $[a,b]$, where $[a,b]=\{ x\in X: a\leq x \leq b \}$ is bounded. 
In what follows, all linearly ordered coarse spaces are suppose to be connected.

 We recall that a subset $Y$ of a  linearly ordered set $(X, \leq)$ is called {\it convex} if $[a,b]\subseteq  Y$ for all  $a,b\in Y$, and observe that $Y$ is convex if and only if there exists $x\in Y$ such that  $[x,y]\subseteq Y$  for each $y\in X$.

 \vspace{7 mm}
{\bf Theorem 1. } {\it For a  coarse space 
$(X, \mathcal{E})$ and a linear  order $\leq$ on $X$, the following  statements are equivalent 

$(i)$  $ \ (X, \mathcal{E}, \leq)$ is linearly ordered;

$(ii)$ $ \  \mathcal{E}$ has a base  
$   \mathcal{E}^\prime$ such  that $E^\prime [x]$ is convex for all $x\in X$, $E^\prime \in \mathcal{E}^\prime$.

\vspace{3 mm}

 Proof. } $(i)$   $\Rightarrow \ (ii) \ $. For each 
 $E\in \mathcal{E}$, we denote 
 $E^\prime= \ \bigcup  \ \{[x,y]: \  x,y\in E\} $.
 Since $E^\prime$ is convex, it suffices to show that 
 $E^\prime\in \mathcal{E}$. 
 Since $\leq$ is compatible with $\mathcal{E}$, there exists $F\in \mathcal{E}$ such that $E\subseteq F$ and if  $x< z$  $\ (z< x)$ and $z\in X\setminus F[x]$ then $y< z$ $ \ (z< y)$ for each $y\in E[x]$. 
 It follows that $E^\prime [x] \subseteq F[x]$  and 
 $E^\prime \in \mathcal{E}$.
 \vspace{3 mm}
 
 $(ii)\Rightarrow (i)$. Given $E \in \mathcal{E}$, we choose $E^\prime \in \mathcal{E}$  such that $E \in E^\prime$ and $E^\prime [x]$ is convex for each $x\in X$.
 If $z\in X\setminus E^\prime [x]$ then either $z< x^\prime$ for each $x^\prime \in E[x]$ or $x^\prime < z$ for each $x^\prime \in E[x]$.
 Hence, $\leq$ is compatible with $\mathcal{E}$.
 $ \ \  \  \Box $
 \vspace{6 mm}

We say that a base $E^\prime$ of  $\mathcal{E}$, satisfying $(ii)$ is {\it locally convex}.

Let $\mathcal{B}$ be a bornology on a set $X$. 
Following \cite{b1}, we say that a coarse structure $\mathcal{E}$ on $X$ is {\it compatible} with 
$\mathcal{B}$ if $\mathcal{B}$ is the bornology of bounded subsets of the coarse space $(X, \mathcal{E})$.

For a linear order $\leq$  on a set $X$, $ \ \  \mathcal{B}_\leq$ denotes the interval bornology on $X$ with the base $\{[a,b]: a,b\in X \}$.
In the following two examples, we describe the smallest 
locally convex coarse structure 
$\downarrow\mathcal{E}_\leq$ and the strongest  locally convex coarse structure $\uparrow\mathcal{E}_\leq$
compatible with $\mathcal{B}_\leq$.

\vspace{6 mm}

{\bf Example 1. } 
 Let $\leq$ be a linear order on a set $X$.
 Then $\downarrow\mathcal{B}_\leq$
is the discrete coarse structure on $X$ defined by the bornology $\mathcal{B}_\leq$.

\vspace{6 mm}

{\bf Example 2. } 
 Let $\leq$ be a linear order on a set $X$, 
 $\mathcal{C}_\leq$ denotes the family of all bounded convex subsets of $(X, \leq)$.
 We consider the family $\Phi$ of all mappings $\varphi: X\rightarrow \mathcal{C}$ such that, for all 
 $a,b\in X$, we have 
 $$\bigcup  \ \{\varphi(x): x\in [a,b] \} \in \mathcal{B}_\leq , \  \  \{ x\in X: \varphi (x) \ \bigcap \  [a,b] \ \neq\varnothing\}\in \mathcal{B}_\leq.$$

 Then the family $\{E_\varphi  :  \varphi\in \Phi \}$,
 where $E_\varphi = \{ (x,y): y\in \varphi(x) \}$, is a base for $\uparrow \mathcal{E}_\leq$.
 
 \vspace{6 mm}
 
 Following  \cite{b7}, we say that a coarse structure $\mathcal{E}$ on $(X, \leq)$ is {\it interval} if there is a base 
 $\mathcal{E}^\prime$ of $\mathcal{E}$ 
 such that,
 for all $E^\prime$, $x\in X$, $E^\prime [x]$ is an interval in 
 $(X, \leq)$.
 Clearly, $\mathcal{E}$ is locally convex and, by Theorem 1, $\leq$ is compatible with $\mathcal{E}$.
 On the other hand, let $(X, \mathcal{E}, \leq)$ be a linearly ordered coarse space. Is $\mathcal{E}$ an interval coarse structure? We give the negative answer to this question. 
 
 The following example also shows  that a subspace of a coarse space with interval base 
may not have  an interval base. 
 
 \vspace{6 mm}

{\bf Example 3. } We denote by $X$ the subset $\bigcup \{(2^n -  1, \  2^n +  1): \  n> 1\}$ of $\mathbb{R}$, put $E_0 = \{(x,y)\in X\times X: \  |x-y| < 2\}$.
$$E_n = \{(x,y)\in X\times X: \  x,y \in (3, \ 2^{n+1})\}, \ n>1,$$ endow $X$ with a coarse structure $\mathcal{E}$ with the base $\{E_n \cup E_0 : n> 1 \}$.
Then $\mathcal{E}$ is locally convex. To  see that $\mathcal{E}$ does not have an interval base, we observe that if $H\subseteq X\times X$, $H[x] $ is an interval for each $x\in X$ and $E_0\subseteq H$ then $2^{n+1} \in H[2^n]$  for each $n>1$. Hence, $H\notin \mathcal{E}$.

\section{ Asymptotic dimension}

Let $(X, \mathcal{E})$ be a coarse space, $E\in \mathcal{E}$. 
A family $\Im$ of subsets of $X$ is called $E$-bounded ($E$-disjoint) if, for each $A\in \Im$, there exists $x\in X$ such that $A \subseteq E[x] \ $
($E[A]\cap B = \emptyset$ for all distinct $A, B \in \Im$).

By the definition [13, Chapter 10], 
$asdim (X, \mathcal{E})\leq n$ if, for each  $E\in \mathcal{E}$, there exist $F\in \mathcal{E}$ and $F$-bounded covering $\mathcal{M}$ of $X$ which can be partitioned $\mathcal{M}= \mathcal{M}_0 \cup \dots \cup \mathcal{M}_n$ so that each family $\mathcal{M}_i$ is $E$-disjoint. 
If there is the minimal natural number $n$ with this property then $asdim (X, \mathcal{E})=n$, otherwise $asdim (X, \mathcal{E})= \infty$.

\vspace{7 mm}
{\bf Theorem 2. } {\it Let $(X, \mathcal{E}, \leq)$ be a linearly ordered coarse space. Then $asdim (X, \mathcal{E})\in  \{0,1\}$.
\vspace{5 mm}

 Proof. } Let $E\in \mathcal{E}$, $E=E^{-1}$
 and $E[x]$ is convex for each $x\in X$, see Theorem 1. 
 We fix $x\in X$, observe that $E^n [x]$ is convex for each $n\in  \mathbb{N}$ and put $E^\omega [x]= \bigcup \{E^n [x]: x\in \mathbb{N}\}$.
 We show that there exist $E^2$-bounded covering $\mathcal{M}(x)$ of $E^\omega [x] $  and a partition
 $\mathcal{M}(x)= \mathcal{M}_0 (x)\bigcup \mathcal{M}_1 (x)$
such that $\mathcal{M}_0 (x), \  \mathcal{M}_1 (x)$ are 
$E$-disjoint.

For each $n\in \mathbb{N}$, we denote 
$R_n = (E^{n+1} [x] \setminus E^{n} [x] )  \  \bigcap \ \{y\in X: x\leq y\} $, 
$ \ \  L_n = (E^{n+1} [x] \setminus E^{n} [x] ) \  \bigcap \ \{y\in X: y< x\} $ and observe that 
$E[R_i]\bigcap L_j =\emptyset$ for all $i,j\in \mathbb{N}$ and 
$E[R_i]\bigcap R_j = \emptyset$, 
$E[L_i]\bigcap L_j = \emptyset \ $ for all $ \ i,j\in \mathbb{N} \ $ such that $ \ \ |i-j|> 1$. 
Clearly, $R_n \ $ and $ \ L_n$
are convex for each $ \ n\in \mathbb{N}$.
If $R_n = \emptyset \ $  ($ \ L_n = \emptyset$ ) then 
$R_i = \emptyset \ $  ($ \ L_i = \emptyset$ ) for every 
$i>n$.

We put
$$\mathcal{M}_0 (x)= \{ E[x], \  R_{2n},  L_{2n}: n\in \mathbb{N} \},$$
$$\mathcal{M}_1 (x)= \{R_{2n-1}, L_{2n-1}: n\in \mathbb{N} \}$$ and note that $\mathcal{M}_0 (x)$, 
$\mathcal{M}_1 (x)$ are $E$-disjoint.

We show that each $R_n$ is $E^2$-bounded, the case $L_n$ is analogous. 
For $n=1$, we have $R_1 \subseteq E^2 [x]$. Let $n>1$ and  $R_n \neq \emptyset$.
We take $y\in R_{n-1}$ such that $E[y]\cap R_n \neq\emptyset$ and show that
$R_n \subseteq E^2 [y]$. 
Given $z\in R_n$, we choose $t\in R_{n-1}$ such that 
$z\in E[t]$.
Since $R_{n-1}$  is convex and $E[y]\cap R_n \neq\emptyset$, we have $t\in E[y]$ so $z\in E^2 [y]$.

To conclude the proof, we choose a subset $Z$ of $X$ such that $\bigcup \{E^\omega [z]: z\in Z \}=X$ and
$E^\omega [z] \bigcap E^\omega [z^\prime] = \emptyset$
for all distinct $z, z^\prime \in Z$. For each $z\in Z$
and 
$E^\omega [z]$, we use above construction to choose 
$\mathcal{M}_0 (z)$ and 
$\mathcal{M}_1 (z)$. We put 
$\mathcal{M}_0 = \bigcup \{ \mathcal{M}_0 (z): z\in Z \}$, 
$\mathcal{M}_1 = \bigcup \{ \mathcal{M}_1 (z): z\in Z \}$.
Then 
$\mathcal{M}_0, \   \mathcal{M}_1$ are $E^2$-bounded, 
  $\mathcal{M}_0, \   \mathcal{M}_1$ are $E$-disjoint.  
  $ \ \  \  \Box $
 \vspace{7 mm}

For every discrete coarse structure $\downarrow \mathcal{E}_\leq$ on a  linearly ordered set
 $(X, \leq)$, we have $asdim (X,  \ \mathcal{E}_\leq, \  \leq) \ = \ 0$.
 
 Let $\leq$ be the natural well ordering on $\omega$.
 Then $\uparrow \mathcal{E}_\leq$ coincides with the universal locally finite coarse structure and, by [5, Theorem 1], 
$asdim (\omega, \ \uparrow \mathcal{E}_\leq, \  \leq) \  = \ 1$.

\section{ Selectors}

Let 
 $( X, \   \mathcal{E}, \ \leq)$ be linearly ordered coarse space, $A\subseteq X$, $E\in \mathcal{E}$.
 We say that $a\in A$ is a {\it right (left)
 $E$-end } of $A$ if $x<a$ ($a<x$) for each $x\in A\setminus E[x]$.
 If $a$ is the maximal (minimal) element of $A$ then $a$ is a right (left) $E$-end for each $E\in \mathcal{E}$.

\vspace{7 mm}
{\bf Example 4. } Let   
 $( X, \   \mathcal{E}, \ \leq)$ be linearly ordered coarse space, metrizable by a metric $d$ on $X$, for metrizability of coarse spaces see [11, Chapter 2].
 We take an arbitrary $\varepsilon \ > \ 0$ and show that every right bounded subset $A$ of $X$ has a right $\varepsilon$-end.
 To this end, we take $a_0\in A$. 
 If $a_0$ is not a right   
 $\varepsilon$-end then we choose $a_1 \in A$ such that 
 $a_0< a_1,  \  \  d(a_0, a_1)>\varepsilon$.
 Repeating this procedure, after finite number of steps, we get a right $\varepsilon$-end $a_n$ of $A$.

 \vspace{7 mm}
{\bf Example 5. } Let   
 $( X, \   \mathcal{E}, \ \leq)$ be a discrete coarse space, defined  by  the interval bornology $\mathcal{B}_\leq$ on $(X, \leq)$.
 Let
 $A\subseteq X$, $ \ B\in \mathcal{B}_\leq$ and let there exists $a\in A$ such that $b<a$ for each $b\in \mathcal{B}$.
 Then $A$ has a right $E_B$-end if and only if $A$ has the  maximal element.
 
 \vspace{7 mm}
{\bf Theorem 3. } {\it Let   
 $( X, \   \mathcal{E}, \ \leq)$
  be a linearly ordered coarse space,
$E\in \mathcal{E}$,
$\Im$ be a family of subsets of $X$. If every subsets 
$A\in \ \Im$ has  a right $E$-end then $\Im$ has a selector. 
  
  \vspace{3 mm}
 
 Proof. } For each $A\in \Im$, we take some right $E$-end $f(A)$ of $A$ and show that the mapping 
 $f: \Im \rightarrow X$
 is macro-uniform.
 
 We take an arbitrary $H\in \mathcal{E}$, $H=H^{-1}$ such that $H[x]$ is convex for each $x\in X$.
 Let $Y, Z\in \Im$, $(Y, Z)\in exp H$ and $f(Y)\leq f(Z)$. We take $y\in Y$ such that $y\in H[f(Z)]$.
 If $y< f(Y)$ then, by the convexity of $H[f(Z)]$ and $f(Y)\leq f(Z)$,  we have $f(Y)\in H[f(Z)]$.
 If $y  \geqslant  f(Y)$ 
 then $y\in E[f(Y)]$. 
 Hence, $H[f(Z)] \ \bigcap \ E[f(Y)] \ \neq \ \emptyset$ and $f(Z) \ \in \ HE[f(Y)]$.
$ \ \  \  \Box $
 \vspace{7 mm}
 
Applying Theorem 3 to Example 4, we conclude that the  family of all right bounded subsets of a lineary ordered metric space has a selector.

\section{ Comments and open questions }

 1. Coarse spaces can be considered as asymptotic 
 counterparts of uniform topological spaces, see
 [11, Chapter 1].
 Selectors and orderings of topological spaces, studied in a plenty of papers, take an important place in {\it Topology}, see surveys  \cite{b2}, \cite{b3}, \cite{b4}, \cite{b12}.
 \vspace{5 mm}
 
2. Example 3 answers negatively
Question 1 from  \cite{b7}, Question 4 was answered negatively in \cite{b8}, Questions 2, 3, 5 from \cite{b7} remain open.

 \vspace{5 mm}

 3. In light of Theorem 3, we ask the following question. 

\vspace{5 mm}

{\bf Question 1. } {\it Does the family of all right bounded subsets of a linearly ordered coarse space have a selector?}
 
 \vspace{5 mm}

4. It is well-knows that every linearly ordered topological space is normal. 
For normality of coarse spaces, see [11, Chapter 4].  

\vspace{5 mm}

{\bf Question 2. } {\it  Is every linearly ordered coarse space normal?}

\vspace{5 mm}

5. We conclude with the following question.  

\vspace{5 mm}

{\bf Question 3. } {\it  Is every linearly ordered coarse space asymorphic to a subspace of a linearly ordered coarse space with an interval base?}

\end{document}